\documentclass[12pt]{amsart}

\usepackage{amsfonts,amsmath,latexsym,amssymb,verbatim,amsbsy,amsthm}

\usepackage{graphicx}
\usepackage{nicefrac}
\usepackage{dsfont}
\usepackage{mathrsfs}
\usepackage{cancel}
\usepackage[normalem]{ulem}
\usepackage{color}

\usepackage[top=1in, bottom=1in, left=1in, right=1in]{geometry}

\usepackage[dvipsnames]{xcolor}
\usepackage{relsize}

\usepackage[colorlinks=true, pdfstartview=FitV, linkcolor=RoyalBlue,citecolor=ForestGreen, urlcolor=blue]{hyperref}

\numberwithin{equation}{section}
\numberwithin{figure}{section}

\renewcommand{\geq}{\geqslant}

\renewcommand{\leq}{\leqslant}

\newcommand{\ds}{\displaystyle} 
\newcommand{\be}{\begin{equation}}
\newcommand{\ee}{\end{equation}}

\theoremstyle{plain}
\newtheorem{THEOREM}{Theorem}[section]

\newtheorem{theorem}[THEOREM]{Theorem}
\newtheorem{corollary}[THEOREM]{Corollary}

\theoremstyle{definition}

\theoremstyle{remark}

\newtheorem{remark}[THEOREM]{Remark}

\newcommand{\myr}[1]{{{#1}}} 

\def \a {\alpha} 
\def \b {\beta}
\def \lam {\lambda}

\def\bm{{\mathbf m}}

\def\rhob_{\rho_{\b}}
\def\rhobp{\rho'_{\b}}
\def\ea{e_{\a}}

\def\fa{f_{\a}}

\def\ubar{\overline{\bu}}

\def\om{\Omega}
\newcommand{\myangle}[1]{(1+{#1})} 

\def \bx {{\mathbf x}}
\def \bxp {{\mathbf x}'}
\def \by {\bxp} 
\def \xp {{x}'}
\def \bu {{\mathbf u}}
\def \bv {{\mathbf v}}
\def \bvp {{\mathbf v}'}

\def\A{{\mathbf A}} 
\def\Etotal{{\mathscr E}} 

\def\delE{\delta\Etotal}
\def\delbu{\delta\bu}

\def\Omegaa{{\mathcal S}_\a}
\def\Omegab{{\mathcal S}_\b}
\def\betaeta{\eta}

\newcommand{\R}{\ensuremath{\mathbb{R}}}   

\def\d{\textnormal{d}}    

\def\e{{\sf{e}}}

\def \hf {\frac{1}{2}}

\def \dx  {\, \d\bx} 
\def\dxp {\, \d\bxp} 
\def\dvp {\, \d\bv'} 
\def \dy  {\dxp} 

\def \dv  {\, \mathrm{d}\bv}

\def \ddt  {\frac{\d}{\d t}} 

\def \Pressure {{\mathbb P}}

\def\pressure{{\raisebox{.4ex}{${}_{\mathbb P}$}}\hspace*{-0.1ex}} 
\def \align {{\mathbf A}}

\newcommand{\bq}{\textbf{q}}

\newcommand{\karray}{k}

\newcommand{\Karray}{K}
\newcommand{\Parray}{\Phi}
\newcommand{\DmK}{\Delta_{\!\!{}_\mathscr{M}}\!\Karray}
\newcommand{\DmP}{\Delta_{\!\!{}_\mathscr{M}}\!\Parray}

\newcommand{\I}{\mathcal{I}}

\begin{document}

\title[Hydrodynamic alignment with pressure II. Multispecies]{Hydrodynamic alignment with pressure II. Multispecies}

\author{Jingcheng Lu}
\address{Department of Mathematics, University of Maryland, College Park}
\email{jlu1@umd.edu}

\author{Eitan Tadmor}
\address{Department of Mathematics and Institute for Physical Sciences \& Technology (IPST), University of Maryland, College Park}
\email{tadmor@umd.edu}

\date{\today}

\subjclass{35Q35, 76N10, 92D25.}

\keywords{swarming, $p$-alignment, multi-species, connectivity, flocking}

\thanks{\textbf{Acknowledgment.} Research was supported by ONR grant N00014-2112773.}

\dedicatory{\bigskip{\large Dedicated to Constantine Dafermos with Friendship and Admiration}}

\begin{abstract}
We study the long-time  hydrodynamic behavior of  systems of multi-species  which arise from agent-based description of alignment dynamics. The interaction between  species is governed by an array of  symmetric communication kernels. 
 We prove that the crowd of different species flock towards the mean velocity if (i) cross-interactions  form a heavy-tailed connected array of kernels, while (ii) self-interactions are governed by kernels with singular heads.  The main new aspect here, is that flocking behavior holds without closure assumption on the specific form of pressure tensors. Specifically, we prove the long-time flocking behavior for connected arrays of multi-species,  with self-interactions governed by entropic pressure laws \cite{Tad2022} and driven by fractional $p$-alignment. In particular, it follows that  such multi-species hydrodynamics  approaches a mono-kinetic description. This  generalizes the mono-kinetic,  ``pressure-less'' study in \cite{HT2021}.
\end{abstract}

\maketitle
\setcounter{tocdepth}{1}
\tableofcontents

\section{Introduction --- alignment dynamics of multi-species}
\subsection{Hydrodynamic description of multi-species}
We study the long-time behavior of the   multi-species hydrodynamics
\begin{subequations}\label{eqs:hydro} 
\begin{equation}\label{eq:hydro}
    \left\{
    \begin{split}
    & \ \partial_{t}\rho_\a+\nabla_\bx\cdot(\rho_\a\bu_\a) = 0,\\
    & \ \partial_{t}(\rho_\a\bu_\a)+\nabla_\bx\cdot(\rho_\a\bu_\a\otimes\bu_\a+\Pressure_\a) = \align_\a(\rho,\bu),
    \end{split}      
    \right. \quad (t,\bx)\in (\R_t, \R^d),
\end{equation}
subject to initial data $\displaystyle (\rho_\a,\bu_\a,\Pressure_\a)_{|_{t=0}} = (\rho_{\alpha0},\bu_{\alpha0},\Pressure_{\alpha0})$.
The different species, tagged by the (possibly infinite) index-set $\a\in {\mathcal I}$, are quantified by their density, $\rho_\a: \R_t\times \R^d \mapsto \R_+$, momentum, $\rho_\a\bu_\a: \R_t\times \R^d \mapsto \R^d$, and pressure tensor, $\Pressure_\a: \R_t\times \R^d \mapsto \R^d\times \R^d$. 
Each species occupies a distinct `patch' of mass $\ds M_\a(t)$ supported on $\Omegaa(t)$,
\[
M_\a(t)=\int \limits_{\Omegaa(t)}\rho_\a(t,\bx)\dx, \qquad \Omegaa(t):=\textnormal{supp} \rho_\a(t,\cdot).
\] 
The dynamics is driven by inter-species
interactions due to \emph{alignment}, dictated by a symmetric array of symmetric communication kernels, $\Phi=\{\phi_{\a\b}(\bx,\bxp)\}$,
\begin{equation}\label{eq:align}
    \begin{split}
        \align_\a(\rho,\bu) &:= \sum_{\b\in {\mathcal I}}\,\int \limits_{\Omegaa(t)}\phi_{\a\b}(\bx,\by)(\bu_{\b}(t,\by)-\bu_{\a}(t,\bx))\rho_{\a}(t,\bx)\rho_{\b}(t,\by)\dy.
    \end{split}
\end{equation} 
 Thus, what distinguishes species $\a$ is the way it communicates with the other species, through symmetric kernels $\phi_{\a\b}, \ \b\in {\mathcal I}$, 
\begin{equation}\label{eq:symmetry}
\phi_{\a\b}(\bx,\by)=\phi_{\a\b}(\by,\bx), \qquad \phi_{\a\b}(\bx,\by)=\phi_{\b\a}(\bx,\by),
\end{equation}
while self-interactions within the same species are governed by $\phi_{\a\a}, \ \a\in{\mathcal I}$.\newline
There is a special role for metric kernels where  communication is dictated by the distance $|\bx-\by|$. In this context we  assume the  existence of  a symmetric array of radially decreasing kernels, $\Karray:=\{k_{\a\b}(r)\}$, such that 
\begin{equation}\label{eq:radial}
\phi_{\a\b}(\bx,\by)\geq k_{\a\b}(|\bx-\by|), \quad k_{\b\a}=k_{\b\a}\geq 0, \qquad 
 \a,\b\in{\mathcal I}.
\end{equation}
\end{subequations}
We use the standard notation $\Phi\succeq \Karray$ to abbreviate \eqref{eq:radial}. This covers the prototypical case of \emph{metric kernels}, $\phi_{\a\b}(\bx,\by)=k_{\a\b}(|\bx-\by|)$, with decreasing intensity of communication as a function of  the distance, e.g., $\phi_{\a\b}(r)=(1+r)^{-\betaeta}$ in \cite{CS2007a}.
In particular, we address general non-decreasing metric kernels, $\phi_{\a\b}(|\cdot|)$, in terms of their \emph{decreasing envelope} $k_{\a\b}(r):=\min\{\phi_{\a\b}(|\bx|) \ | \ |\bx|\leq r\}$. The variety of different classes of communication kernels reflect  large literature on collective dynamics which arises in different disciplines, \cite{Aok1982,VCBCS1995,CF2003,CS2007a,CDMBC2007, Bal2008,CFL2009,Ka2011,MT2011,GWBL2012, MCEB2015,JJ2015,LZTM2019,MLK2019, ST2020b,ST2021}.\newline
The different species are viewed as moving `patches' of  different crowds with mass and momentum which interact according  to the alignment protocol \eqref{eqs:hydro}. We make the following three assumptions about these `patches'.
We assume that  the density of species inside their `patch' remains strictly bounded away from vacuum, 
\begin{equation}\tag{H1}\label{eq:vacuum}
\min_{\bx\in\Omegaa(t)}\rho_\a(t,\bx)\geq \rho_->0, \qquad \forall \a\in{\mathcal I}.
\end{equation}
Further, we assume that 
\begin{equation}\tag{H2}\label{eq:H2} 
\Omegaa(t) \ \textnormal{have smooth boundary satisfying a Lipschitz or a cone condition}, \quad \forall\a\in {\mathcal I}.
\end{equation}
Finally, we assume that the   boundary of each patch
forms a contact discontinuity, governed by Neumann boundary conditions 
\begin{equation}\tag{H3}\label{eq:BCs}
\bu_\a\cdot  {{\mathbf n}_\a}_{|\partial\Omegaa}=0, \quad \Pressure_\a{{\mathbf n}_\a}_{|\partial\Omegaa}=0 \ \ \textnormal{and} \ \ \bq_\a\cdot{{\mathbf n}_\a}_{|\partial\Omegaa}=0, \quad \forall\a\in {\mathcal I}.
\end{equation}
In particular, it follows  that there is no flux of mass  for each species: integration of  \eqref{eq:hydro}${}_1$  implies the mass of each species is conserved
\begin{equation}\label{eq:conserve-mass}
M_\a(t) = M_{\a 0}, \qquad M_\a(t)=\int \limits_{\Omegaa}\rho_\a(t,\bx)\dx.
\end{equation}
In particular, the total mass is also conserved
$\ds M:=\sum \limits_\a M_\a(t)=\sum \limits_\a M_{\a0}$.
In contrast, the momentum of each species need not necessarily conserved due to the cross alignment terms between different species on the right of  \eqref{eq:hydro}${}_2$. Instead, the symmetry of $\phi_{\a\b}(\cdot,\cdot)$ implies that the \emph{total momentum} is conserved\footnote{Here and below we abbreviate
$\square':=\square(t,\by)$.}
\[
\begin{split}
\ddt \sum_\a &\int \limits_{\Omegaa}\rho_\a\bu_\a\dx \\
 & =-\int \limits_{\partial\Omegaa} \left(\bu_\a\cdot{\mathbf n}_\a\bu_\a +\Pressure_\a{\mathbf n}_\a\right){\d}S + \sum_{\a,\b}\ \ \iint \limits_{(\bx,\by)\in \Omegaa\times\Omegab} \phi_{\a\b}(\bx,\by)(\bu'_\b-\bu_\a)\rho_\a\rho_\b' \dx\dy=0,
\end{split}
\] 
and hence
\begin{equation}\label{eq:conserve-momentum}
\bm:=\sum_\a \bm_\a(t) = \sum_\a \bm_{\a 0}, \qquad \bm_\a(t):=\int \limits_{\Omegaa} \rho_\a(t,\bx)\bu_\a(t,\bx)\dx.
\end{equation}

\subsection{The class of entropic pressure laws}
The multi-species system \eqref{eq:hydro} requires a closure for the pressure tensors $\Pressure_\a(t,\bx), \a\in {\mathcal I}$. 
 In this context, we  recall the notion of \emph{entropic pressure} \cite{Tad2022}. 
We refer to  $\Pressure_\a$ as an \emph{entropic pressure} tensor associated with species $\a$ in \eqref{eq:hydro} if its non-negative trace $\rho_\a e_\a(t,\bx):=\hf \text{trace}(\Pressure_\a(t,\bx))\geq 0$ satisfies
\begin{equation}\label{eq:meso-pressure}
    \partial_{t}(\rho_{\a}e_{\a})+\nabla_{\bx}\cdot(\rho_{\a}e_{\a}\bu_{\a}+\bq_{\a})+\textnormal{trace}(\Pressure_\a\nabla\bu_\a) \leq - 2 \sum_\b \int \limits_{\Omegab} \phi_{\a\b}(\bx,\by) e_\a\rho_\a\rho'_{\b}\dy.
\end{equation}
Here $\bq_\a(t,\bx)$ is an $C^1$-flux.\newline 
The motivation for \eqref{eq:meso-pressure} stems from 
the large-crowd dynamics of the agent-based model proposed in \cite{HT2021},
in which different species, each of which consists of $N_\a$ agents with position/velocity  $\displaystyle (\bx^{\a}_{i}(t),\bv^{\a}_{i}(t)): \R_{+}\mapsto \R^{d}\times\R^{d}$, are driven by the Cucker-Smale alignment \cite{CS2007a}
\begin{equation}
\label{eq:multi CS}
    \left\{\begin{array}{c}
    \begin{split}
        \ddt\bx^{\a}_{i}(t) &= \bv^{\a}_{i}(t),   \\
        \ddt \bv^{\a}_{i}(t) &= \sum_{\b\in {\mathcal I}}\frac{1}{N_{\b}}\sum_{j = 1}^{N_{\b}}\phi_{\a\b}(\bx^{\b}_{j}(t),\bx^{\a}_{i}(t))(\bv^{\b}_{j}(t)-\bv^{\a}_{i}(t)),
    \end{split}      
    \end{array}\right. \qquad i=1,2,\ldots N_\a, 
\end{equation}
The passage from the agent-based to the hydrodynamic description goes through an intermediate kinetic description which is realized by the  empirical distribution, 
$\ds     f_{\a}(t,\bx,\bv) := \frac{1}{N_{\a}}\sum_{i=1}^{N_{\a}}\delta_{\bx_{i}^{\a}(t)}\otimes\delta_{\bv_{i}^{\a}(t)}$.
Indeed, the hydrodynamic description  \eqref{eq:hydro} is recovered in terms of  the first-two limiting moments of $\{f_\a\}$ which are assumed to exist,  
\[
\rho_\a(t,\bx)= \lim_{N_{\a}\rightarrow\infty}\int \limits_{\R^d} f_{\a}(t,\bx,\bv)\dv, \quad 
\rho_\a\bu_\a(t,\bx) = \lim_{N_{\a}\rightarrow\infty}\int \limits_{\R^d} \bv f_{\a}(t,\bx,\bv)\dv.
\] 
This process of  large-crowd limit as $N_\a\rightarrow \infty$   recovers \eqref{eq:hydro} with pressure, $\Pressure_\a$, given by the  second-order moments
\begin{equation}\label{eq:pressure}
        \Pressure_{\a}(t,\bx) = \lim_{N_{\a}\rightarrow\infty}\int \limits_{\R^d} (\bv-\bu_\a)\otimes (\bv-\bu_\a)f_\a(t,\bx,\bv)\dv.
\end{equation}
The formal derivation is outlined in appendix   \ref{sec:hydro-des} and it follows the different derivations  with different level of rigor in  case of single species
\cite{HT2008,CFTV2010,CCR2011,FK2019,NP2021,Shv2021,NS2022}.
The kinetic description of the pressure in terms of second-order, rank-one moments in \eqref{eq:pressure} leads to the notion of \emph{internal energy} which quantifies  microscopic  fluctuations around the bulk velocity $\bu_\a$,
\[
\rho_\a e_\a = \hf \text{trace}(\Pressure_\a)= \lim_{N_{\a}\rightarrow\infty} \int \limits_{\R^d} \hf |\bv-\bu_\a|^2f_\a(t,\bx,\bv)\dv.
\]
This kinetic description of internal energy  formally yields the equality
\[
\partial_{t}(\rho_{\a}e_{\a})+\nabla_{\bx}\cdot(\rho_{\a}e_{\a}\bu_{\a}+\bq_{\a})+\textnormal{trace}(\Pressure_\a\nabla\bu_\a) = - 2 \sum_\b \int \limits_{\Omegab} \phi_{\a\b}(\bx,\by)e_\a\rho_\a \rho'_{\b}\dy,
\]
with heat flux $\ds \bq_{\a}:=\lim_{N_{\a}\rightarrow\infty} \int \limits_{\R^d} \hf |\bv-\bu_\a|^2(\bv-\bu_\a)f_\a(t,\bx,\bv)\dv$.
Thus, we arrive at the special case of equality in \eqref{eq:meso-pressure}.
In particular, it covers the  ``pressure-less'' case --- the special case of \emph{mono-kinetic closure} 
\[
f_\a(t,\bx,\bv) \stackrel{N_\a\rightarrow \infty}{\longrightarrow}\rho_\a(t,\bx)\delta(\bv-\bu_\a(t,\bx)),
\]
 which is realized in terms of zero pressure, $\Pressure_\a=0$, 
\begin{equation}\label{eq:mono-hydro}
    \left\{\begin{array}{c}
    \begin{split}
    & \partial_{t}\rho_\a+\nabla_\bx\cdot(\rho_\a\bu_\a) = 0,\\
    & \partial_{t}(\rho_\a\bu_\a)+\nabla_\bx\cdot(\rho_\a\bu_\a\otimes\bu_\a) = \align_\a(\rho,\bu),
    \end{split}      
    \end{array}\right. \quad (t,\bx)\in (\R_t, \R^d),
\end{equation}
Most of the literature on swarming hydrodynamics of single species \emph{assumes}  mono-kinetic closure. The corresponding ``pressure-less'' multi-species  hydrodynamics was studied  in \cite{HT2021}.
The definition of pressure in terms of the entropy \underline{inequality} \eqref{eq:meso-pressure} is not concerned, however,  with the precise details of internal energy, as it lacks a reference to the specific closure with respect to a preferred state of thermal equilibrium. In fact, \eqref{eq:meso-pressure} applies to a large class of tensors beyond those which are  realizable as second-order moments.

\subsection{Energy dissipation in entropic alignment}\label{sec:energy-dissipation}
 A main consequence of the notion of entropic pressure  is to secure the dissipativity of the total energy $\displaystyle  E_\a:=\frac{|\bu_\a|^2}{2} +e_\a$.
Indeed, manipulating the mass and momentum equations we find 
 \begin{equation}\label{eq:energy-kin}
 \begin{split}
 \partial_t \big(\frac{\rho_\a}{2}|\bu_\a|^2\big) + &\nabla_\bx\cdot \big(\frac{ \rho_\a}{2}|\bu_\a|^2\bu_\a +\Pressure_\a\bu_\a\big) - \textnormal{trace} \big(\Pressure_\a \nabla\bu_\a\big) \\
  & = 
 - \sum_\b \int\limits_{\Omegab}\phi_{\a\b}(\bx,\by)(|\bu_\a|^2-\bu_\a\cdot\bu'_\b)\rho_\a\rho'_\b\dy.
 \end{split}
 \end{equation}
Adding the entropic description of the pressure postulated in \eqref{eq:meso-pressure} we end up with, 
\begin{equation}\label{eq:energy-inequality}
 \begin{split}
 \partial_t \big(\rho_\a E_\a\big) +  &\nabla_\bx\cdot \Big(\rho_\a E_\a\bu_\a +\Pressure_\a\bu_\a +\bq_\a\Big)  \\
  & \leq 
 - \sum_\b \int\limits_{\Omegab}\phi_{\a\b}(\bx,\by)\big(|\bu_\a|^2-\bu_\a\cdot\bu'_\b+2 e_\a\big)\rho_\a\rho'_\b\dy.
 \end{split}
 \end{equation}
Thus,  the role of entropic pressure is to complement the energy balance  \eqref{eq:energy-kin} in forming an entropy inequality \eqref{eq:energy-inequality}, which augments  the system of hyperbolic balance laws \eqref{eqs:hydro}; we refer to the authoritative book of  \cite{Daf2016}.
This implies dissipativity of the total energy. Indeed,  by the zero Neumann boundary conditions  assumed in \eqref{eq:BCs}, it follows that
\begin{equation}\label{eq:energy-ineq}
\begin{split}
\ddt \sum_\a \int \limits_{\Omegaa(t)} & \rho_\a E_\a\dx  \\
 & \leq
-\sum_{\a,\b}\,\iint \limits_{\Omegaa(t)\times\Omegab(t)} \phi_{\a\b}(\bx,\bxp)
\big(|\bu_\a|^2-\bu_\a\cdot\bu'_\b+2 e_\a\big)\rho_\a\rho'_\b\dx\dy\\
 & = -\hf \sum_{\a,\b}\,\iint \limits_{\Omegaa(t)\times\Omegab(t)} \phi_{\a\b}(\bx,\bxp)
\big(|\bu'_\b-\bu_\a|^2+2e_\a+2e'_\b\big)\rho_\a\rho'_\b\dx\dy.
\end{split}
\end{equation}
For further discussion we refer to \cite[\S1]{Tad2022}.

\section{Swarming and long-time flocking behavior}
  We discuss the behavior of a large crowd,  possibly infinite number  of species, $\{\rho_\a,\bu_\a,\Pressure_\a\}$, each of which consists of a large crowd of agents, $\{\bx_i^\a,\bv_i^\a\}$ in \eqref{eq:multi CS}. A crowd of species  (or agents) is viewed as a \emph{swarm} when it is driven by collective dynamics which coordinates  its species (or agents) to aggregate together with emergence of   large-scale formations. 
  In the present context  of  dynamics  governed by alignment \eqref{eqs:hydro}, we  are concerned with the long-time \emph{flocking behavior} of the multi-species system \eqref{eqs:hydro}. Flocking refers to  the emergence of coherent structure with limiting velocities  $\bu_\a^\infty$ such that
\[
\bu_\a(t,\bx)-\bu_\a^\infty(t,\bx) \stackrel{t\rightarrow \infty}{\longrightarrow}0,
\]
 with the corresponding limiting densities, $\rho_\a^\infty=\rho_\infty(\bx-\bu_\a^\infty t)$.   Since we ignore attraction, repulsion or external forcing,  the limiting behavior of pure alignment should be particularly simple --- the different species governed by  \eqref{eq:hydro} can only approach the same time-invariant mean velocity 
 \[
\bu_\a(t,\cdot) \stackrel{t\rightarrow \infty}{\longrightarrow}\ubar, \qquad \ubar:=\frac{\sum_\a \bm_\a}{\sum_\a M_\a},
\]
 with  a limiting density  carried out as a traveling wave $\rho_\a^\infty(\bx-\ubar t)$. Ideally, we seek uniform convergence. In the present context multi-species with pressure, we have no access to uniform bounds on the velocities. Instead a more relaxed notion  of $L^2_\rho$-convergence becomes accessible by studying \emph{energy fluctuations}, 
\begin{equation}\label{eq:eflock}
\delE(t):=\sum_{\a} \ \int \limits_{\Omegaa}\left\{\hf |\bu_{\a}(t,\bx)-\ubar|^2+e_\a(t,\bx)\right\}\rho_\a(t,\bx)\dx.
\end{equation}
Observing that 
\begin{equation}\label{eq:delEeqE}
\delE(t) =\sum_\a \int \limits_{\Omegaa} \rho_\a E_\a\dx -\langle \ubar, \sum_\a {\bm_\a}\rangle 
+ \hf |\ubar|^2 M =\sum_\a \int \limits_{\Omegaa} \rho_\a E_\a\dx - \hf |\ubar|^2 M,
\end{equation}
we conclude that  energy fluctuations decay at the same  rate  as the total energy in \eqref{eq:energy-ineq}
\begin{equation}\label{eq:efluc-decay}
\ddt \delE(t) \leq -\hf \sum_{\a,\b}\ \,\iint \limits_{\Omegaa(t)\times\Omegab(t)} \phi_{\a\b}(\bx,\bxp)
\big(|\bu'_\b-\bu_\a|^2+2e_\a+2e'_\b\big)\rho_\a\rho'_\b\dx\dy.
\end{equation}
Our flocking results will be quantified in terms of  the decay of  energy/energy fluctuations, which in turn implies the decay of both --- the macroscopic velocity fluctuations around the \emph{mean velocity} $\ubar$, and the microscopic (kinetic) fluctuations of the different species around their bulk velocities, $\ds \rho_\a e_\a = \lim_{N_{\a}\rightarrow\infty} \int \limits_{\R^d} \hf |\bv-\bu_\a|^2f_\a(t,\bx,\bv)\dv$. \newline
A second component of flocking behavior requires that alignment is strong enough to keep the dynamics  contained in a finite ball, forming the `flock' 
 \[
  D(t):=\sum_\a D_\a(t)\leq D_+<\infty, \qquad D_\a(t):=\max_{\bx,\bxp \in {\mathcal S}_\a(t)}|\bx-\bxp|.
  \]
In practice we may need to  address  a relaxed notion of flocking which allows a slow time growth, $D(t) \leq C_D\myangle{t}^{\gamma}$ with some fixed $\gamma>0$.
  
\subsection{Statement of main results}
  The multi-species alignment dynamics  \eqref{eq:hydro} is dictated by  the array of communication kernels $\Phi=\{\phi_{\a\b}\}$. Our flocking results require $\Phi$ to form a connected array. To this end, it will suffice to consider the radial lower bounds $ \phi_{\a\b}(\bx,\bxp)\geq \karray_{\a\b}(|\bx-\bxp|$
  assumed in \eqref{eq:radial}.
 The array  $\Karray=\{\karray_{\a\b}\}$ is viewed as the adjacency matrix of  a \emph{weighted graph}, with a {weighted graph Laplacian}, $\DmK(r)$, \cite{HT2021}
\[
(\DmK(r))_{\a\b} := \left\{\begin{array}{lr}
    \ds -\karray_{\a\b}(r)\sqrt{M_{\a}M_{\b}}, &\hspace{0.5em} \alpha \neq \beta,\\ \\
    \ds \sum_{\gamma\neq\a} \karray_{\a\gamma}(r)M_{\gamma}, & \hspace{0.5em} \alpha = \beta.
    \end{array}\right.
\]
Algebraic connectivity is quantified in terms of the \emph{spectral gap},  $\lambda_2(\DmK)$, \cite{Fie1973,Fie1989},
\begin{equation}\label{eq:gap}
 \lambda_2(\DmK):= M\min_{\mathbf{y}}\left\{\sum_\a\sum_{\b\neq \a}\karray_{\a\b}|y_\a-y_\b|^2M_{\a}M_{\b} \ \Big| \ \sum_\a\sum_{\b\neq \a}|y_\a-y_\b|^2M_{\a}M_{\b}=1\right\}.
 \end{equation}
The graph associated with $\Karray$ is connected if and only if $\lambda_2(\DmK)>0$. Since the spectral gap  is a non-decreasing  function of  the non-negative entries, \cite[\S3]{HT2021},  $\lambda_2(\DmK)>0$ also  controls the connectivity of the  communication kernels, $\lambda_2(\DmP)\geq \lambda_2(\DmK)>0$. 
  Our flocking results  require \emph{heavy-tailed connectivity} in the sense that $\lambda_2(\DmK(r))$ has slow enough decay in a manner made precise in  the theorem  \ref{thm:main1} below.

\smallskip\noindent
{\bf Notations}. Below, we use $C_K$ and $C_D$ to denote constants which characterize the heavy-tailed behavior of $K$ and the dispersion of diameter $D(t)$. We let $C_R$  
denote a constant, with different values in different contexts, depending of a spatial scale $R$, as well as on the other fixed parameters \myr{on} the problem $\betaeta,\gamma,...$. Finally, we let $C_1,C_2,\ldots$ \myr{denote} related parameters which arise from computations with these constants. 

\begin{theorem}\label{thm:main1}
Consider the multi-species system \eqref{eqs:hydro} with two or more species.
Let $(\rho_\a,\bu_\a,\Pressure_\a)$ be a non-vacuous strong solution\footnote{That is, \eqref{eq:vacuum} holds for $(\rho_{\a}(t,\cdot),\bu_{\a}(t,\cdot),\Pressure_{\a}(t,\cdot))\in(L^{\infty}\cap L^{1}_{+}(\R^{d}))\times W^{1,\infty}(\R^{d})\times W^{1,\infty}(\R^{d})$, $\alpha\in {\mathcal I}$} of  \eqref{eqs:hydro}, subject to compactly supported initial data $(\rho_{{\a}0},\bu_{{\a}0},\Pressure_{{\a}0})$ with $D(0)<\infty$, and boundary conditions \eqref{eq:BCs}. 
 Assume that $\Karray(r)$ has heavy-tailed connectivity of order $\betaeta\geq0$,  namely --- there exist $C_K, R>0$ such that
\begin{equation}\label{eq:connect}
   \lambda_2(\DmK(r)) \geq C_K (1+r)^{-\betaeta},   \qquad  r\geq R.
\end{equation}
Moreover, assume that the crowd disperses  at the rate of order $\gamma\geq0$, namely --- there exists $C_D>0$ such that for all $t\geq 0$,
\begin{equation}\label{eq:sizeofDt}
D(t) \leq C_D(1+t)^\gamma, \qquad \gamma\geq 0.
\end{equation}
 If the heavy-tail condition holds in the sense that
 \be\label{eq:heavy_tail}
 \betaeta\gamma < 1,
 \ee
  then there is a large time flocking behavior  with fractional exponential decay rate 
\begin{equation}\label{eq:flocking_decay}
\delE(t) 
  \leq C_R\, exp\big\{-C_\zeta t^{(1-\betaeta\gamma)}\big\}\delE(0), \qquad C_\zeta:=2\zeta C_KC_D^{-\betaeta}, \quad \zeta:=1-\frac{\max_\a M_\a}{M}>0. 
\end{equation}
\end{theorem}
Theorem \ref{thm:main1} extends the mono-kinetic, pressure-less case \cite[Theorem 4.1]{HT2021}. It applies to general entropic pressure laws \eqref{eq:pressure}, and general symmetric communication protocol satisfying \eqref{eq:connect}: the kernels $\phi_{\a\b}$ need not be metric nor upper-bounded. At the same time, it extends   the heavy-tail condition for flocking of a single species asserted in \cite[corollary 3.3]{Tad2022}.\newline
The decay estimate \eqref {eq:flocking_decay} reflects  a competition between the possible dispersion of the crowd as its diameter  $D(t)$ may grow in time, and the decay rate in  the strength of communication strength, $\lam_2\big(\DmK(r)\big)$, as the `edge of the crowd' may grow with  $r$. Theorem \ref{thm:main1} tells us that if their composition has  a  non-integrable tail so that
\[
\int \limits_0^t \lam_2\big(\DmK(D(\tau)\big){\d}\tau \leq  C_R e^{-C_\zeta 1t^{(1-\betaeta\gamma)}} \stackrel{t\rightarrow \infty}{\longrightarrow}0, \qquad \betaeta\gamma<1,
\]
 then the  different species flock towards  the mean velocity $\displaystyle \overline{\bu}_\infty$
\[
\sum_{\a\in {\mathcal I}}\int|\bu_{\a}(t,\bx)-\overline{\bu}_{\infty}|^{2}\rho_{\a}(t,\bx)\dx \lesssim \e^{-C_\zeta1 t^{(1-\betaeta\gamma)}}\delE(0); 
\]
 moreover,  there is a (fractional) exponential decay of internal fluctuations, 
\[
\sum_\a \int \|\Pressure_\a(t,\bx)\|^2\dx= \sum_{\a\in {\mathcal I}}\int|\bv-\bu_{\a}(t,\bx)|^{2}f_\a(t,\bx,\bv)\dv\dx \lesssim e^{-C_\zeta t^{(1-\betaeta\gamma)}}\delE(0).
\]

\subsection{The example of ``pressure-less'' equations} A key aspect of theorem \ref{thm:main1} is a dispersion bound which controls the spatial diameter, $D(t)\lesssim \myangle{t}^\gamma$. As a prototypical example we consider the  mono-kinetic ``pressure-less'' closure, $\Pressure_\a=0$, \cite{HT2021}. In this case, the alignment dynamics \eqref{eq:hydro}${}_2$ decouples into scalar transport equations for the components of $\bu_\a=(u_\a^1, \ldots, u_\a^d)$, 
\[
\partial_t u^i_\a + \bu_\a\cdot \nabla_\bx u_\a^i= 
\sum_\beta \int \limits_{\Omegab(t)}\phi_{\a\b}(\bx,\bxp)\big(u_\b^i(t,\bxp)-u_\a^i(t,\bx)\big)\rho'_\b\dxp.
\]
Assume $\betaeta<1$, then a maximum principle of the scalar velocity components eventually leads to the uniform bound $D(t)\leq D_+<\infty$, i.e., theorem \ref{thm:main1} applies with $\gamma=0$, leading to exponential decay
$\ds \delE(t)   \leq C_Re^{-C_\zeta t}\delE(0)$. In fact, there is  exponential decay  of velocity fluctuations in the \emph{uniform} norm \cite[step \#3 in the proof of theorem 1.1]{HT2021}
\[
\max_\a \max_{\bx\in\Omegaa(t)}|\bu_\a(t,\bx)-\ubar| \leq  C_Re^{-C_\zeta t}\max_\a \max_{\bx\in\Omegaa(t)}|\bu_{\a0}(\bx)-\ubar|.
\]
We conclude that flocking of  pressure-less  dynamics is dictated for any heavy-tailed connectivity of order $\betaeta<1$, \eqref{eq:connect}
\[
\lambda_2(\DmK(r)) \geq C_K (1+r)^{-\betaeta},   \qquad  \betaeta<1, \quad r\geq R.
\]

\section{Self-interactions based on fractional $p$-alignment}
We now turn our attention to the main aspect of this work --- multi-species alignment with pressure. In this case, one does not have  access to pointwise bounds on the velocities $\bu_\a$, which in turn imply the desired pointwise bound on the diameters, $D_\a(t)$, propagating with these velocities.  Instead, we follow the single-species  arguments of \cite[\S6]{Tad2022}, in order to secure  direct bounds the dispersion of $D_\a(t)$. To this end, observe that the heavy-tailed  flocking scenario in theorem \ref{thm:main1}   is quantified in terms of the spectral gap  \eqref{eq:gap} which is \emph{independent of}  self interactions, $\{\phi_{\a\a}\}$.  The desired dispersion bound will be obtained when we consider enhanced self-interactions; specifically  ---  we consider self-interactions based on \emph{singular communication kernels}, 
$\ds \phi_{\a\a}(\bx,\bxp)=|\bx-\bxp|^{d+2sp}, \ 0<s<1, \ p\geq 1$.
Such kernels greatly emphasize the alignment with immediate neighbors over far away neighbors, leading to 
\begin{equation}\label{eq:phydro}
\left\{\begin{split}
\ds \partial_t&(\rho_\a \bu_\a) + \nabla_\bx\cdot (\rho_\a \bu_\a\otimes \bu_\a+\Pressure_\a) = \\ \\
&= \int \limits_{\Omegaa}\frac{|\bu'_\a-\bu_\a|^{2p-2}(\bu'_\a-\bu_\a)}{|\bxp-\bx|^{d+2sp}}\rho_\a\rho'_\a{\d}\bxp
+\sum_{\b\neq \a}\, \int \limits_{\Omegab}\phi_{\a\b}(\bx,\bxp)(\bu'_\b-\bu_\a)\rho_\a\rho'_\b{\d}\bxp.
\end{split}\right.
\end{equation}
Self-interactions in this case amount to \emph{weighted fractional $2p$-Laplacians}; more precisely, the first integrand on the right of  \eqref{eq:phydro} is the subdifferential of the weighted Gagliardo fractional energy (suppressing the time dependence) \cite{DPV2012}
\[
{\mathcal J}_{2p,s}(\bu_\a)= \iint \frac{|\bu_\a(\bxp)-\bu_\a(\bx)|^{2p}}{|\bxp-\bx|^{d+2sp}}\rho_\a \rho'_\a\dx\dxp, \qquad 0<s<1, \ p\geq 1.
\]
Interactions  based on $p$-alignment, $p>1$ in the context of a single species were  introduced  in  \cite{HKK2014,CCH2014} and further developed in \cite{Tad2022}. We note that \eqref{eq:phydro} corresponds to the multi-species agent-based description with self interactions based on $p$-alignment \begin{equation}\label{eq:pCS}
\ddt \bv^{\a}_{i} = \frac{1}{N_{\a}}\sum_{j = 1}^{N_{\a}}\phi_{\a\a}(\bx^{\a}_{j},\bx^{\a}_{i})|\bv^\a_j-\bv^\a_i|^{2p-2}\big(\bv^{\a}_{j}-\bv^{\a}_{i}\big) + \sum_{\b\neq \a}\frac{1}{N_{\b}}\sum_{j = 1}^{N_{\b}}\phi_{\a\b}(\bx^{\b}_{j},\bx^{\a}_{i})(\bv^{\b}_{j}-\bv^{\a}_{i}),
\end{equation}
with singular kernels $\phi_{\a\a}(\bx,\bxp)=|\bx-\bx|^{-(d+2sp)}$.
The passage from \eqref{eq:pCS} to \eqref{eq:phydro} can be justified only in the case of bounded (or at least integrable) $\phi_{\a\a}$. and remains formal in the singular case.\newline
We close this section noting that since the flocking bound in \eqref{eq:flocking_decay} is independent of self-interactions, the main theorem \ref{thm:main1} still applies to  the case  of self-interactions based on fractional $p$-alignment in \eqref{eq:phydro}.
\subsection{Energy dissipation in entropic $p$-alignment}\label{sec:energy-pdissipation}
The notion of `entropic pressure' in \eqref{eq:meso-pressure}  requires an adjustment for $p$-alignment. Following \cite[remark 6.1]{Tad2022}, 
we refer to  $\Pressure_\a$ as an \emph{entropic pressure} tensor associated with species $\a$ in \eqref{eq:phydro} with $C^1$ `heat-flux'  $\bq_\a$, if its non-negative trace $\rho_\a e_\a:=\hf \text{trace}(\Pressure_\a)\geq 0$ satisfies
\begin{equation}\label{eq:meso-ppressure}
\begin{split}
    \partial_{t}&(\rho_{\a}e_{\a})+\nabla_{\bx}\cdot(\rho_{\a}e_{\a}\bu_{\a}+\bq_{\a})+\textnormal{trace}(\Pressure_\a\nabla\bu_\a) \\
     & \leq -  \hf  D_\a^{-(d+2sp)}(t)\int  \limits_{\Omegaa} \big((2e_\a)^p+(2e'_\a)^p\big)\rho_\a\rho'_\a \dy
     -2  \sum_{\b\neq \a}\, \int \limits_{\Omegab} \phi_{\a\b}(\bx,\by)e_\a\rho_\a\rhobp \dy.
     \end{split}
\end{equation}
The self-interaction terms in \eqref{eq:meso-ppressure} and \eqref{eq:meso-pressure} for `pure' alignment $p=1$  end with the same energy dissipation statement. Specifically, manipulating \eqref{eq:hydro}${}_1$ and \eqref{eq:phydro} yields, corresponding to \eqref{eq:energy-kin},
\[
 \begin{split}
 \partial_t \big(\frac{\rho_\a}{2}|\bu_\a|^2\big) + \nabla_\bx\cdot \big(\frac{ \rho_\a}{2}|\bu_\a|^2\bu_\a +&\Pressure_\a\bu_\a\big) - \textnormal{trace} \big(\Pressure_\a \nabla\bu_\a\big)   \\
    = & - \int \limits_{\Omegaa} \frac{|\bu'_\a-\bu_\a|^{2p-2}\bu_\a\cdot(\bu_\a-\bu'_\a)}{|\bxp-\bx|^{d+2ps}}\dx \\
 &  - \sum_{\b\neq \a} \,\int\limits_{\Omegab}\phi_{\a\b}(\bx,\by)(|\bu_\a|^2-\bu_\a\cdot\bu'_\b)\rho_\a\rho'_\b\dy.
 \end{split}
\]
Adding the entropic description of the pressure postulated in \eqref{eq:meso-ppressure} we find, arguing along the lines \eqref{eq:efluc-decay}
\[
\begin{split}
\ddt \delE(t) & \leq - \hf\sum_\a \, \iint\limits_{\Omegaa(t)\times \Omegaa(t)}\Big(\frac{|\bu'_\a-\bu_\a|^{2p}}{|\bxp-\bx|^{d+2sp}}+ D_\a^{d+2sp}(t)\big((2e_\a)^p+(2e'_\a)^p\big)\Big)\rho_\a\rho'_\a\dx\dy \\
  &  \ \ \ -\hf \sum_{\b \neq \a} \ \, \iint \limits_{\Omegaa(t)\times\Omegab(t)} \phi_{\a\b}(\bx,\bxp)
\big(|\bu'_\b-\bu_\a|^2+2e_\a+2e'_\b\big)\rho_\a\rho'_\b\dx\dy.
\end{split}
\]
In particular, ignoring the negative  contributions coming from internal energy and  from cross interactions terms, yields
\begin{equation}\label{eq:energy-pdissipation}
 \ddt \delE(t)      \leq  
  - \hf\sum_\a \, \iint\limits_{\Omegaa(t)\times \Omegaa(t)}\frac{|\bu'_\a-\bu_\a|^{2p}}{|\bxp-\bx|^{d+2sp}}\rho_\a\rho'_\a\dx\dy.
 \end{equation}
Thus, the contribution coming from self interactions based on singular $p$-alignment  imply that the velocities $\bu_\a$ are bounded in the (homogeneous) Sobolev spaces $\dot{W}^{s,2p}(\Omegaa)$. Specifically, taking into account the non-vacuous bound assumed in  \eqref{eq:vacuum} then integration of \eqref{eq:energy-pdissipation}  yields
\begin{equation}\label{eq:enstrophy}
\begin{split}
\int \limits_0^t \sum_\a \|&\bu_\a(\tau,\bx)\|^{2p}_{\dot{W}^{s,2p}(\Omegaa)}{\d}\tau  \\
 &   \leq 
  C^2_\rho \int \limits_0^t \sum_\a \iint \limits_{\Omegaa(t)\times \Omegaa(t)} \frac{|\bu_\a(t,\bxp)-u_\a(t,\bx)|^{2p}}{|\bxp-\bx|^{d+2sp}}\rho_\a\rho'_\a\dx\dy{\d}\tau\\
   & \qquad \quad\leq   C^2_\rho C^2_0, \qquad \qquad C^2_0:=2\sum_\a \int \limits_{\Omegaa(0)} \rho_{\a0} E_{\a0} \dx, \quad C_\rho:=\frac{1}{\rho_-}.
   \end{split}
\end{equation}
\subsection{Multi-species with entropic pressure and fractional $p$-alignment }
The enstrophy bound \eqref{eq:enstrophy} implies a dispersion bound sought in \eqref{eq:sizeofDt}. We follow the argument in \cite{Tad2022}.
The mass propagation by \eqref{eq:hydro}${}_1$ implies
\[
\ddt D_\a(t) \leq \delbu_\a(t), \qquad \delbu_\a(t):=\max_{\bx,\bxp \in \Omegaa(t)}|\bu(t,\bx)-\bu(t,\bxp)|.
\] 
Gagliardo-Nirenberg inequality implies for all $\nicefrac{d}{2p} < s<1$ there holds,\footnote{It is here that we use the assumed smoothness of the boundaries of ${\Omegaa}$ in \eqref{eq:H2}.} \cite{DPV2012,MRR2013}
\[
|\bu(t,\bx)-\bu(t,\bxp)| \leq C_s\|\bu_\a(t,\cdot)\|_{\dot{W}^{s,2p}(\Omegaa(t))}|\bx-\bxp|^{s-\theta}, \qquad \theta:=\frac{d}{2p}<s <1,
\]
and hence $\ds \ddt D_\a(t) \leq C_s\|\bu_\a(t,\cdot)\|_{\dot{W}^{s,2p}(\Omegaa)}D_\a^{s-\theta}(t)$.
It follows,   that
\begin{equation}\label{eq:revisit}
\ddt \sum_\a D^{1+\theta-s}_\a(t) \leq C'_s \sum_\a \|\bu_\a(t,\cdot)\|_{\dot{W}^{s,2p}(\Omegaa)}, \qquad C'_s= (1+\theta-s)C_s.
\end{equation}
\myr{Now, since $1+\theta-s <1$, then $D^{1+\theta-s}_\a \leq \sum_\a D^{1+\theta-s}_\a$,  and integration of \eqref{eq:revisit} yields,}
\[
\begin{split}
D^{1+\theta-s}(t) &\leq \sum_\a D_{\a 0}^{1+\theta-s} + \Big(\int \limits_0^t \sum_\a\|\bu_\a(\tau,\cdot)\|^{2p}_{\dot{W}^{s,2p}(\Omegaa)}{\d}\tau\Big)^{\tfrac{1}{2p}}
\Big(\int \limits_0^t 1{\d}\tau\Big)^{\tfrac{1}{(2p)'}} \\
 & \leq \sum_\a D_{\a 0}^{1+\theta-s} + C_s'(C_\rho C_0)^{\tfrac{1}{p}} t^{\tfrac{1}{(2p)'}}.
 \end{split}
 \]
We conclude that  multi-species  crowd driven by self-interaction of $p$-alignment dynamics, 
\eqref{eq:phydro} can be dispersed at a rate no faster than
\begin{equation}\label{eq:Dabeta}
D(t) \leq C_D(1+t)^{\gamma_p}, \qquad \gamma_p= \frac{2p-1}{2p(1+\theta-s)}, \quad \theta=\frac{d}{2p}<s <1.
\end{equation}
This bound can be improved: in appendix \ref{sec:dispersion} we use a bootstrap argument to show a slower rate of order
\[
D(t) \leq C'_D(1+t)^{\myr{\gamma_*}}, \qquad \gamma_*= \frac{\myr{2p-1}}{2p(1+\theta-s)+\betaeta}, \quad \theta=\frac{d}{2p}<s <1
\]
Theorem \ref{thm:main1} applies, leading to flocking behavior  of order $\lesssim exp\{-t^{1-\betaeta\gamma_*}\}$ which we summarize in the following.
  
\begin{theorem}\label{thm:main2}
Let $(\rho_\a,\bu_\a,\Pressure_\a)$ be a non-vacuous strong solution of  \eqref{eq:hydro}${}_1$,\eqref{eq:phydro} satisfying \eqref{eq:vacuum}--\eqref{eq:BCs}, with cross interactions, $\Phi(\bx,\bxp) \succeq \Karray(|\bx-\bxp|)$, and self-interactions based on  $p$-alignment of order  $p>\tfrac{d}{2}$.
 Assume that $\Karray(r)$ has tail connectivity of order $\betaeta\geq0$, \eqref{eq:connect} 
 \[
   \lambda_2(\DmK(r)) \geq C_K (1+r)^{-\betaeta},   \qquad  r\geq R.
\]
 If the heavy-tail condition \eqref{eq:heavy_tail} holds,
 \[
 \betaeta\gamma_p < 1, \qquad \gamma_p:=\frac{2p-1}{2p(1+\theta-s)}, \quad \theta=\frac{d}{2p}<s <1,
 \]
  then there is a large time flocking behavior  with fractional exponential decay rate 
\begin{equation}\label{eq:flocking_decayB}
\delE(t) 
  \leq C_R\, exp\big\{-C_\zeta t^\mu\big\}\delE(0), \quad \mu= {\frac{2p(1+\theta-s)-\myr{2(p-1)}\betaeta}{2p(1+\theta-s)+\betaeta}}>0, 
\end{equation}
with constant $C_\zeta=2\zeta C_K(C'_D)^{-\betaeta}$. 
\end{theorem}
\begin{remark}[{\bf Lack of exponential decay bound}] We leave open the question of  a uniform dispersion bound, $D(t)\leq D_+<\infty$ corresponding to $\gamma_*=0$, which in turn would imply the \emph{exponential} decay
$\ds \delE(t) \leq C_Re^{-C_\zeta t}$.  This will require an improved bootstrap argument in appendix \ref{sec:dispersion}, along the lines of \cite[Appendix D]{Tad2022}. 
\end{remark}
\subsection{Multi-species with entropic pressure in one-dimension} The methodology leading to theorem \ref{thm:main2} consists of two main parts: (i) a $\betaeta$-\emph{tailed} array of cross-interactions which enforce flocking  of multi-species dynamics; and (ii)   self-interactions based on $p$-alignment with singular \emph{head} which dictate the dispersion rate $\gamma_p$. Observe that this line of argument requires $\nicefrac{d}{2p}<1$, and therefore flocking of self-interactions based on `pure' alignment, $p=1$, is restricted to the $d=1$-case,
\begin{equation}\label{eq:1Dhydro}
\left\{\begin{split}
\ds \partial_t&(\rho_\a u_\a) + \partial_x(\rho_\a u^2_\a+\pressure_\a) = \\ \\
&= \int \limits_{\Omegaa}\frac{(u'_\a-u_\a)}{|\xp-x|^{1+2s}}\rho_\a\rho'_\a{\d}\xp
+\sum_{\b\neq \a}\, \int \limits_{\Omegab}\phi_{\a\b}(x,\xp)(u'_\b-u_\a)\rho_\a\rho'_\b{\d}\xp,
\end{split}\right. \  (t,x)\in (\R_t, \R)
\end{equation}
with scalar entropic pressures, $\pressure_\a$, satisfying (assuming  no heat flux $q_\a=0$),
\begin{equation}\label{eq:scalar-p}
\partial_t  \pressure_\a + \partial_x(\pressure_\a u)+2\pressure_\a\partial_x  u \leq -2\pressure_\a D_\a^{1+2s}(t)M.
\end{equation}
Theorem \ref{thm:main2} applies  with $\ds \gamma_1=\frac{1}{3-2s}$.
\begin{corollary}[{\bf Multi-species in one-dimension}]\label{cor:main2}
Consider the one-dimensional multi-species system \eqref{eq:hydro}${}_1$,\eqref{eq:1Dhydro} with entropic pressure \eqref{eq:scalar-p} and satisfying \eqref{eq:vacuum},\eqref{eq:BCs}.  
 If the heavy-tail connectivity condition holds 
 \[
 \betaeta+2s <3 , \quad \hf<s<1,
 \]
   then there is a large time flocking behavior  with fractional exponential rate 
\begin{equation}\label{eq:1Dflocking_decay}
\delE(t) 
  \leq exp\Big\{-2C_\zeta t^\mu\Big\} \delE(0), \qquad \mu={\frac{3-2s}{3-2s+\betaeta}}>0. 
\end{equation}
\end{corollary} 
Singular  interactions of a  single species in one dimension  with mono-kinetic closure were extensively studied in \cite{ST2017a,ST2017b,DKRT2018,ST2018a,ST2020b} and we refer to the review \cite{MMPZ2019} and the additional references therein. Corollary \ref{cor:main2}  extends  these flocking results to multi-species in one dimension with  entropic pressure laws.
Going beyond the one-dimensional corollary \ref{cor:main2},  clarifies the motivation for our discussion of  self-interactions based on fractional $p$-alignment, $p>1$, which extend the discussion  to higher dimensions.

\appendix
\section{From agent-based to hydrodynamic description}\label{sec:hydro-des} 
We begin with the  derivation of the multi-species hydrodynamic description \eqref{eq:hydro} from the agent-based dynamics \eqref{eq:multi CS}.

The large crowd dynamics of of the different species can be encoded in terms of their empirical distribution $\displaystyle f_\a(t,\bx,\bv):= \frac{1}{N_\alpha}\sum \limits_{i=1}^{N_\a} \delta_{\bx^\a_i(t)}(\bx)\otimes \delta_{\bv^\a_i(t)}(\bv)$, which are governed by the kinetic Valsov equation in state variables   $(t,\bx,\bv)\in \R_+\times\om\times \R^d$, 
 e.g., \cite{HT2021},
\begin{equation}\label{eq:Q}
\begin{split}
\partial_t f_\a +\bv\cdot\nabla_\bx f_\a + \nabla_\bv\cdot Q_\a(f_\a,{\mathcal F})=0,\qquad {\mathcal F}=\{f_\b\},
 \end{split}
\end{equation}
where different species are interconnected through pairwise  communication protocol on the right (we abbreviate $f_\a=f_\a(t,\bx,\bv), f'_\b=f_\b(t,\bxp,\bvp)$ and likewise $\square=\square(t,\bx), \ \square'=\square(t,\bxp)$ etc.)  
\[
\begin{split}
 Q_\a(f_\a,{\mathcal F})  :=\sum_\b \int\limits_{\Omegab}\phi_{\a\b}(\bx,\bxp)(\bvp-\bv)f_\a f'_\b \dvp\dxp.
 \end{split}
\]
The large crowd dynamics of  $f_\a$'s is captured by their first two moments which we assume to exist --- the density $\displaystyle \rho_\a(t,\bx):=\lim_{N_\a\rightarrow \infty} \int \limits_{\R^d} f_\a(t,\bx,\bv)\dv$, and   the momentum, $\displaystyle  \rho_\a \bu_\a(t,\bx):=\lim_{N_\a\rightarrow \infty} \int \limits_{\R^d} \bv f_\a(t,\bx,\bv)\dv$. Integration of \eqref{eq:Q} yields the mass equation \eqref{eq:hydro}${}_1$ 
\[
\partial_{t}\rho_\a+\nabla_\bx\cdot(\rho_\a\bu_\a) = 0.
\]
The first $\bv$-moment of \eqref{eq:Q} yields 
\begin{equation}\label{eq:ki-momentum}
\partial_t \int \limits_{\R^d} \bv f_\a\dv  =-\nabla_\bx\cdot \int \limits_{\R^d} \bv \bv^\top f_\a\dv  +   \int \limits_{\R^d} Q_\a(f_\a, {\mathcal F})\dv.
  \end{equation} 
  For the first term on the right of \eqref{eq:ki-momentum} 
  $\bv \bv^\top\equiv  -\bu_\a\bu_\a^\top+(\bv\bu^\top_\a+\bu_\a\bv^\top)+ (\bv-\bu_\a) (\bv-\bu_\a)^\top$, where the first two moments of $f_\a$ add up to  $\bu_\a (\rho\bu)_\a^\top= \rho_\a\bu_\a \otimes \bu_\a$, and the third yields the pressure tensor \eqref{eq:pressure}, 
 \[
 \int \limits_{\R^d} \bv \bv^\top f_\a\dv =
 \rho_\a\bu_\a \otimes \bu_\a +  \Pressure_\a, \qquad \Pressure_\a=\int \limits_{\R^d} (\bv-\bu_\a)(\bv-\bu_\a)^\top f_\a\dv;
 \]
 the second term on the right of \eqref{eq:ki-momentum} yields
 \[
 \int \limits_{\R^d} Q_\a(f_\a,{\mathcal F})\dv = \sum_\b \int \limits_{\Omegab(t)} \phi_{\a\b}(\bx,\bxp)\left((\rho\bu)'_\b\rho_\a-(\rho\bu)_\a\rho'_\b\right)\dxp =\A_\a(\rho,\bu),
\]
and we recover the momentum equation \eqref{eq:hydro}${}_2$
\[
 \partial_{t}(\rho_\a\bu_\a)+\nabla_\bx\cdot(\rho_\a\bu_\a\otimes\bu_\a+\Pressure_\a) = \align_\a(\rho,\bu).
\]
Observe that the system \eqref{eqs:hydro} is not a purely hydrodynamic description  since  the pressure    in \eqref{eq:pressure} still requires a  \emph{closure} of the $\bv$-dependent second-order moments of $f_N$. This is our point of departure from the flocking analysis in \cite{HT2021}:  the   hydrodynamic description of alignment in \eqref{eqs:hydro}, is  left open. Following \cite{Tad2022}, we will trace the  decay of energy fluctuations, showing that  it   applies to general entropic  pressure stress tensors \eqref{eq:pressure}. 

\subsection{Energy balance}
We derive the energy balance  as preparation for studying the long-time behavior of hydrodynamics \eqref{eqs:hydro}. The total energy is given by the second moment which is assumed to exist
\begin{equation*}
    \rho_{\a}E_{\a}(t,\bx) = \lim_{N_{\a}\rightarrow\infty}\int \limits_{\R^d}\frac{|\bv|^{2}}{2}\fa(\bx,\bv,t)\dv,
\end{equation*}
is decomposed into kinetic and internal energy corresponding to the decomposition $\displaystyle \hf|\bv|^{2} = \hf|\bu_\a|^{2} + \hf|\bv-\bu_\a|^{2}+\bu_\a\cdot(\bv-\bu_\a)$ and noticing that $\ds \int_{\R^d} (\bv-\bu_\a)f_\a\dv=0$,  
\begin{equation*}
\rho_\a E_\a = \frac{\rho_\a}{2}|\bu_\a|^{2}+\rho_\a e_\a, \qquad \rho_\a e_\a:= \hf\int \limits_{\R^d} |\bv-\bu_\a|^2f_\a\dv.
\end{equation*}
The balance of energy can be obtained by integrating (\ref{eq:Q}) against $\displaystyle \frac{|\bv|^{2}}{2}$, obtaining
\[
\partial_t (\rho_\a E_\a) + I_\a = II_\a.
\]
Here $I_\a$  is the transport-based term which we express as
\[
\begin{split}
 I_\a & = \int \limits_{\R^d} \frac{|\bv|^{2}}{2}(\bv\cdot\nabla_{\bx}f_\a)\dv  \\
& \equiv \nabla_{\bx}\cdot \int \limits_{\R^d} \frac{|\bv|^{2}}{2}\bu_{\a}\fa\dv+\nabla_{\bx}\cdot\int\frac{|\bv|^{2}}{2}(\bv-\bu_{\a})\fa\dv 
\\
 & \equiv \nabla_{\bx}\cdot \left(\,\int \limits_{\R^d} \frac{|\bv|^{2}}{2}\fa\dv\right)\bu_{\a}
  +\nabla_{\bx}\cdot\int \limits_{\R^d} \left[\frac{|\bu_\a|^2}{2} + (\bv-\bu_{\a})\cdot \bu_{\a}+\frac{|\bv-\bu_\a|^2}{2}\right](\bv-\bu_\a)f_\a\dv\\
  & =  \nabla_{\bx}\cdot (\rho_{\a}E_{\a}\bu_{\a}) +
  \nabla_\bx \cdot\int \limits_{\R^d} (\bv-\bu_\a)(\bv-\bu_\a)^\top\bu_\a f_\a\dv +
  \nabla_{\bx}\cdot\int \limits_{\R^d} \frac{|\bv-\bu_{\a}|^{2}}{2}(\bv-\bu_{\a})f_\a\dv\\
        & = \nabla_{\bx}\cdot(\rho_{\a}E_{\a}\bu_\a+\Pressure_{\a}\bu_{\a}+\bq_{\a}), 
\end{split}
\]
involving the pressure tensor
$\Pressure_\a$, \eqref{eq:pressure}, and a heat-flux vector, $\bq_\a$,
\begin{equation}\label{eq:heat-flux}
\bq_{\a}:=\hf\int \limits_{\R^d} |\bv-\bu_{\a}|^{2}(\bv-\bu_{\a})\fa\dv,
\end{equation}
and  $II_\a$ is the   alignment-based enstrophy term   given by
 \[
 \begin{split}
  II_\a  & = \int \limits_{\R^d} \bv\cdot Q_\a(f_\a, {\mathcal F})\dv =   -\sum_{\b\in\I}\, \int \limits_{\Omegab(t)}\phi_{\a\b}(\bx,\bx') \Big(\,\iint \limits_{\R^d\times \R^d }\bv\cdot(\bv-\bvp)f_\a f'_\b \dv\dvp\Big)\dy\\
   & =  -\sum_{\b\in\I} \ \int \limits_{\Omegab(t)} \phi_{\a\b}(\bx,\by)(2E_{\a}-\bu_{\a}\cdot\bu'_{\b})\rho_{\a}\rho'_{\b}\dy.
   \end{split}
 \]
Combining,  we formally end up with the energy balance
\begin{equation}\label{eq:energy-balance}
\begin{split}
    \partial_{t}(\rho_{\a}E_{\a})&+\nabla_{\bx}\cdot(\rho_{\a}E_{\a}\bu_{\a}+\Pressure_{\a}\bu_{\a}+\bq_{\a}) \\
     & = -\sum_{\b\in\I} \ \int \limits_{\Omegab(t)} \phi_{\a\b}(\bx,\by)(|\bu_{\a}|^2-\bu_{\a}\cdot\bu'_{\b}+\ea)\rho_{\a}\rho'_{\b}\dy.
     \end{split}
\end{equation}
Thus, the energy equality  which arises form a kinetic description,  is viewed here as a special case of the inequality \eqref{eq:energy-ineq} associated with the general class of entropic pressures,
\begin{equation}\label{eq:energy-ineqb} 
\begin{split}
    \partial_{t}(\rho_{\a}E_{\a}) &+\nabla_{\bx}\cdot(\rho_{\a}E_{\a}\bu_{\a}+\Pressure_{\a}\bu_{\a}+\bq_{\a}) \\
     & = -\sum_{\b\in\I} \ \int \limits_{\Omegab} \phi_{\a\b}(\bx,\by)(|\bu_{\a}|^2-\bu_{\a}\cdot\bu'_{\b}+\ea)\rho_{\a}\rho'_{\b}\dy.
    \end{split}
\end{equation}
\subsection{Energy fluctuations}
Integrating \eqref{eq:energy-ineqb} and summing over 
$\displaystyle \a\in\I$ we find 
 \begin{equation}\label{eq:etfluc}
\begin{split}
\frac{\d}{{\d}t}\sum_{\a\in\I} \ &\int\limits_{\Omegaa(t)}\rho_{\a}E_{\a}\dx\\
&  \leq   -\sum_{\a}\int \limits_{\partial\Omegaa(t)}\big(\rho_{\a} E_{\a}\bu_{\a}\cdot{\mathbf n}_{\a} +  \Pressure_{\a}\bu_{\a} \cdot {\mathbf n}_{\a} +\bq_{\a}\cdot{\mathbf n}_{\a}\big) \d S\\ 
& \ \ \  -\sum_{\a,\b\in\I}\ \ \iint\limits_{\Omegaa(t) \times \Omegab(t)}\phi_{\a\b}(\bx,\bx')(|\bu_{\a}|^2-\bu_{\a}\cdot\bu'_{\b}+\ea)\rho_{\a}\rho'_{\b}\dx\dy\\
    & =  -\sum_{\a,\b\in\I}\ \ \iint \limits_{\Omegaa(t) \times \Omegab(t)}\phi_{\a\b}(\bx,\bx')\left(\hf |\bu_{\a}-\bu'_{\b}|^2+e_{\a}+e'_{\b}\right)\rho_{\a}\rho'_{\b}\dx\dy.
\end{split}
\end{equation}
The boundary conditions assumed in \eqref{eq:BCs} imply there is no energy flux and hence  the boundary integrals on the right vanish\footnote{In fact, here one can consider a larger class  of \emph{energy dissipative} boundary condition}, 
while the symmetrization assumed in \eqref{eq:symmetry} implies, upon change of variables $(\a,\bx)\leftrightarrow(\b,\bxp)$, that the second term  admits the symmetric form of the integrals on the right.  

The inequality \eqref{eq:etfluc} quantifies the energy dissipation in terms of the negative  total enstrophy on the right. This is better expressed in an equivalent symmetric form, in terms of the \emph{energy fluctuations}
\begin{equation}\label{eq:efluc}
\delE(t)=\frac{1}{2M}\sum_{\a,\b\in {\mathcal I}} \ \ \iint \limits_{\Omegaa(t) \times \Omegab(t)} \left\{\hf |\bu_{\a}(t,\bx)-\bu_{\b}(t,\bx')|^2+e_\a(t,\bx)+e_\b(t,\by)\right\}\rho_{\a}\rho'_\b\dx\dy.
\end{equation}
\begin{remark}\label{rem:eremark}
Observe that the definition of energy fluctuation in \eqref{eq:efluc} coincides with the one we had in \eqref{eq:eflock}. Indeed, since the total mass in \eqref{eq:conserve-mass} and momentum in \eqref{eq:conserve-momentum}, are conserved in time,  then  the mean velocity  is invariant in time, 
\[
\ubar=\frac{\sum_\a \bm_\a}{\sum_\a M_\a} = 
\frac{\sum_\a \bm_{\a 0}}{\sum_\a M_{\a 0}},
\]
and the macroscopic   portion of the energy fluctuations \eqref{eq:efluc} can be expressed as fluctuations around that mean velocity, 
\[
\begin{split}
&\sum_{\a,\b}\,\iint  \limits_{\Omegaa(t) \times \Omegab(t)}  \left(\hf|\bu_\a(t,\bx)-\bu_\b(t,\bxp)|^2\right)\rho_\a\rho'_\b\dx\dy \\
 &\sum_{\a,\b}\, \iint  \limits_{\Omegaa(t) \times \Omegab(t)}  \left(\hf|\bu_\a(t,\bx)- \ubar|^2 + (\bu_\a-\ubar)\cdot(\ubar-\bu'_\b)
+\hf|\bu_\b(t,\bxp)-\ubar|^2\right)\rho_\a\rho'_\b\dx\dy\\
& =
M\sum_\a\int \limits_{\Omegaa(t)}  |\bu_\a(t,\bx)-\ubar|^2\rho_\a(t,\bx)\dx. 
\end{split}
\]
Hence,  the energy fluctuations \eqref{eq:efluc} coincides with its equivalent definition  \eqref{eq:eflock} stated in theorem \ref{thm:main1}
\[
\delE(t)=\sum_\a \int \limits_{\Omegaa(t)}|\bu_\a(t,\bx)-\ubar|^2\rho_\a(t,\bx)\dx.
\]
\end{remark}
Noting that 
\[
\delE(t)= \sum_{\a}\int\left\{\hf \rho_{\a}|\bu_{\a}|^2+\rho_\a e_\a\right\}\dx-\frac{1}{2M}\Big|\sum_\a \bm_\a\Big|^2, \quad \bm_\a:=\int \bu_\a\rho_\a\dx,
\]
with a total mass, $M:=\sum_\a M_\a$, and total momentum, $\sum_\a \bm_\a$, which are  conserved in time  we end up with the symmetric version of the dissipation statement \eqref{eq:etfluc}, expressed in terms of energy fluctuations,
\begin{equation}\label{eq:esymmetric}
\begin{split}
\ddt \sum_{\a,\b\in {\mathcal I}}\ \ &\iint  \limits_{\Omegaa(t)\times \Omegab(t)}\left\{\hf |\bu_{\a}-\bu'_{\b}|^2+e_\a+e'_\b\right\}  \rho_{\a}\rho'_\b\dx\dy \\
 & =  \frac{\d}{{\d}t}\sum_{\a\in\I}\ \int\limits_{\Omegaa(t)}\rho_{\a}E_{\a}\dx \\
 &\leq 
- \sum_{\a,\b\in {\mathcal I}}\ \ \iint \limits_{\Omegaa(t)\times \Omegab(t)}\phi_{\a\b}(\bx,\by)\left\{\hf |\bu_{\a}-\bu'_{\b}|^2+e_\a+e'_\b\right\}\rho_{\a}\rho'_\b\dx\dy.
\end{split}
\end{equation}
 
\section{Flocking of strong solutions --- proof of theorem \ref{thm:main1}}
Recall that the metric kernels $k_{\a\b}(r)$ are assumed to decrease with the distance $r$ and hence \eqref{eq:esymmetric} implies
\[
\begin{split}
\frac{\d}{{\d}t}\delE(t)  & \leq -\sum_{\a,\b\in\I} \ \ \iint \limits_{\Omegaa(t)\times \Omegab(t)}\hspace*{-0.2cm}k_{\a\b}(|\bx-\bx'|)\left(\hf |\bu_{\a}-\bu'_{\b}|^2+e_{\a}+e'_{\b}\right)\rho_{\a}\rho'_{\b}\dx\dy\\
& \leq -\sum_{\a,\b\in\I}k_{\a\b}(D(t))\hspace*{-0.2cm} \iint \limits_{\Omegaa(t)\times \Omegab(t)} \left(\hf |\bu_{\a}-\bu'_{\b}|^2+e_{\a}+e'_{\b}\right)\rho_{\a}\rho'_{\b}\dx\dy
\end{split}
\]
For the first term on the right we use the weighted Poincare inequality \cite[Lemma 3.2]{HT2021}
\[
\begin{split}
\sum_{\a,\b\in\I}&k_{\a\b}(D(t))\hspace*{-0.2cm} \iint \limits_{\Omegaa(t)\times \Omegab(t)} \hf |\bu_{\a}-\bu'_{\b}|^2\rho_{\a}\rho'_{\b}\dx\dy\\
& \geq \frac{\zeta}{M} \lambda_2(\DmK(D(t))\sum_{\a,\b\in\I}\ \ \iint \limits_{\Omegaa(t)\times \Omegab(t)} \hf |\bu_{\a}-\bu'_{\b}|^2\rho_{\a}\rho'_{\b}\dx\dy, \quad \zeta=1-\frac{\max_\a M_\a}{M}.
\end{split}
\]
For the remaining terms
\[
\begin{split}
\sum_{\a,\b\in\I}k_{\a\b}(D(t))&\hspace*{-0.2cm} \iint \limits_{\Omegaa(t)\times \Omegab(t)} \left(e_{\a}+e'_{\b}\right)\rho_{\a}\rho'_{\b}\dx\dy \\
& = \sum_\b deg_\b(\Karray) \sum_\a \int \limits_{\Omegaa(t)} e_{\a}\rho_\a \dx + \sum_\a deg_\a(\Karray) \sum_\b \int \limits_{\Omegab(t)}e'_{\b}\rho'_\b \dy,
\end{split}
\]
where $deg_\a$ --- the degree of connectivity of species $\a$, has the lower bound, \cite[eq. (3.10)]{HT2021}
\[
deg_\a(\Karray):= \sum_{\gamma \neq \a} k_{\a\gamma}(D(t)) M_\gamma
\geq \zeta \lambda_2(\DmK(D(t)).
\]
We end up with
\[
\begin{split}
\sum_{\a,\b\in\I}k_{\a\b}(D(t)) &\iint \limits_{\Omegaa(t)\times \Omegab(t)} \left(e_{\a}+e'_{\b}\right)\rho_{\a}\rho'_{\b}\dx\dy \\
& \geq \frac{\zeta}{M}\lambda_2(\DmK(D(t))\sum_{\a,\b\in\I}\ \ \iint \limits_{\Omegaa(t)\times \Omegab(t)} \left(e_{\a}+e'_{\b}\right)\rho_{\a}\rho'_{\b}\dx\dy.
\end{split}
\]
Adding the last two inequalities we conclude the dissipation statement of energy fluctuations
\begin{equation}\label{eq:ms-efluc}
\begin{split}
\frac{\d}{{\d}t}&\delE(t) \\
 & \leq -2\zeta \lambda_2(\DmK(D(t)) \times
\frac{1}{2M}\sum_{\a,\b\in\I} \ \ \iint \limits_{\Omegaa(t)\times \Omegab(t)} \left(\hf |\bu_{\a}-\bu'_{\b}|^2+e_{\a}+e'_{\b}\right)\rho_{\a}\rho'_{\b}\dx\dy \\
& \ \ \ -2\zeta \lambda_2(\DmK(D(t)) \times \delE(t).
\end{split}
\end{equation}
We now turn to address the main flocking bound in \eqref{eq:flocking_decay}.\newline
\begin{proof}[Proof\nopunct] of theorem  \ref{thm:main1}.
Integrate the decay of energy fluctuations  \eqref{eq:ms-efluc},
\[
\ddt \delE(t) \leq -2\zeta\lambda_2(\DmK(D(t)) \times \delE(t),
\]
combined with the assumed bounds \eqref{eq:connect} and \eqref{eq:sizeofDt},
   \[
   \lambda_2(\DmK(D(t)) \geq C_K \big(1+C_D(1+t)^\gamma\big)^{-\betaeta},    \qquad  r\geq R.
   \]
   implies the desired bound
  \[ 
   \delE(t) \leq C_R\, exp\big\{-2\zeta C_KC_D^{-\betaeta} t^{1-\betaeta\gamma}\big\}\delE(0), \qquad \betaeta\gamma<1,
   \]
   with a constant $C_R>0$.
   \end{proof}

\section{An improved dispersion bound}\label{sec:dispersion}
 Assume that we secured the dispersion bound
  $\ds   D(t) \leq C_D(1+t)^\gamma$.
   Then, theorem \ref{thm:main1} applies, leading to flocking behavior  with fractional exponential decay rate which we rewrite as
\begin{equation}\label{eq:flocking_decayA}
\delE(t) 
  \leq C_R\frac{1}{\chi(t)}\delE(0), \qquad \chi(t):= exp\big\{C_\zeta t^{(1-\betaeta\gamma)}\big\}. 
  \end{equation}
    This bound, which was shown to hold with $ \gamma=\gamma_p=\frac{2p-1}{2p(1+\theta-s)}$,  can be improved. To this end, we rewrite \eqref{eq:energy-pdissipation},\eqref{eq:flocking_decayA} in the form
  \[
  \begin{split}
  \ddt \chi(t)\delE(t) &\leq -\frac{\rho_-^2}{2}\chi(t)\sum_\a  \|\bu_\a(t,\cdot)\|^{2p}_{\dot{W}^{s,2p}(\Omegaa)} + C_R\frac{\dot{\chi}(t) }{\chi(t)}\delE(0)\\
  &\leq -\frac{\rho_-^2}{2}\chi(t)\sum_\a  \|\bu_\a(t,\cdot)\|^{2p}_{\dot{W}^{s,2p}(\Omegaa)} + C_2 t^{-\betaeta\gamma}\delE(0), \quad C_2=2C_\zeta(1-\betaeta\gamma)C_R.
  \end{split}
  \]
  This implies
  \[
  \int \limits_0^t \chi(\tau)\sum_\a  \|\bu_\a(\tau,\cdot)\|^{2p}_{\dot{W}^{s,2p}(\Omegaa)}{\d}\tau \leq 2C_\rho^2\delE(0) + C_3t^{1-\betaeta\gamma}\delE(0), \qquad C_3=2C_\rho^2C_2\frac{1}{1-\betaeta\gamma}.
  \]
  Now we revisit \eqref{eq:revisit} with the last weighted bound, obtaining
  \[
  \begin{split}
  D^{1+\theta-s}(t) &\leq \sum_\a D^{1+\theta-s}_{\a 0}
  +  \Big(\int \limits_0^t \sum_\a\chi(t) \|\bu_\a(\tau,\cdot)\|^{2p}_{\dot{W}^{s,2p}(\Omegaa)}{\d}\tau\Big)^{\tfrac{1}{2p}}
\Big(\int \limits_0^t \Big(\frac{1}{\chi^{1/2p}(\tau)}\Big)^{(2p)'}{\d}\tau\Big)^{\tfrac{1}{(2p)'}} \\
 & \leq \sum_\a D_{\a 0}^{1+\theta-s} + C_4 t^{\tfrac{1-\betaeta\gamma}{2p}}, \qquad C_4= \big(C_3\delE(0)\big)^{\tfrac{1}{2p}}\int \limits_0^\infty
 \chi^{-\tfrac{1}{2p-1}}(\tau){\d}\tau.
 \end{split}
 \]
 Thus, we end up with the improved dispersion bound
  \[
    D(t) \leq C'_D(1+t)^{\gamma'}, \qquad \gamma'=\frac{1-\betaeta\gamma}{2p(1+\theta-s)}.
  \]
  This argument can be repeated: since   $\eta\gamma_p <1$ then $\ds \frac{\betaeta}{2p(1+\theta-s)} = \frac{\betaeta\gamma_p}{2p-1}<  1$ and   hence the iterations $\gamma\mapsto \gamma'$ converge to $\ds \myr{\gamma_\infty}=\frac{1}{2p(1+\theta-s)+\betaeta}$. In particular, \myr{since $2p-1>1$ then} after finitely many iterations we  \myr{will reach}   the improved rate, $\myr{\gamma_*=(2p-1)\gamma_\infty}$, 
  \[
  D(t) \leq C'_D(1+t)^{\gamma_*}, \qquad \gamma_*={\frac{\myr{2p-1}}{2p(1+\theta-s)+\betaeta}}.
  \]

\end{document}